\documentclass[A4,10pt]{amsart}
\usepackage{amsmath,amssymb,amsrefs,amsfonts,eucal,yfonts,amscd,comment,mathtools}

\newtheorem{theorem}{Theorem}[section]
\newtheorem{proposition}[theorem]{Proposition}
\newtheorem{lemma}[theorem]{Lemma}
\newtheorem{corollary}[theorem]{Corollary}

\theoremstyle{definition}

\newtheorem{definition}[theorem]{Definition}

\newcommand{\Real}{\mathbb R}                                   
\newcommand{\Natural}{{\mathbb N}}                              
\newcommand{\1}{\boldsymbol{1}}
\newcommand{\0}{\boldsymbol{0}}

\DeclareMathOperator{\ftp}{\overline{\otimes}}
\begin{document}
\title{Tensor products of Archimedean Riesz spaces: a representational approach}

\author{A.W. Wickstead}
\address{Mathematical Sciences Research Centre, Queen's University
Belfast, Belfast BT7 1NN, Northern Ireland.} \email{A.Wickstead@qub.ac.uk}
\begin{abstract}In an earler paper Buskes and the author pointed out that given two Archimedean Riesz spaces $E$ and $F$ it is relatively simple to construct, from their Ogasawara-Maeda representations, a third Archimedean Riesz space $G$ and a bi-injective Riesz bimorphism $\phi:E\times F\to G$. They further pointed out that the Riesz subspace of $G$ geneated by $\phi(E\times F)$ is isomorphic to the Archimedean Riesz space tensor product $E\ftp F$ constructed by Fremlin. The proof there relies on the properties of $E\ftp F$ proved by Fremlin. In this paper we show that this approach actually gives a simple way to construct and establish the properties of $E\ftp F$ and that any choice of $G$ and $\phi$ yield isomorphic objects.  
\end{abstract}
\maketitle

 \section{Introduction}

In 1972 Fremlin, \cite{F}, constructed a tensor product $E\ftp F$ of Archimedean Riesz spaces $E$ and $F$ satisfying the four properties
\begin{enumerate}
\item $E\ftp F$ is an Archimedean Riesz space.
\item $E\otimes F$ is a vector subspace of $E\ftp F$.
\item The Riesz subspace of $E\ftp F$ generated by $E\otimes F$ is the whole of $E\ftp F$.
\item If $G$ is an Archimedean Riesz space and $\phi:E\times F\to G$ is a Riesz bimorphism then there is a unique Riesz homomorphism $T:E\ftp F\to G$ such that $T(x\otimes y)=\phi(x,y)$ for $x\in E$ and $y\in F$.
\end{enumerate}
He also established that this object is essentially unique and gave two fundamental density results. Since then Schaefer, \cite{S}, has given a slightly simplified construction. Both Fremlin and Schaefer's approaches make some use of representations of Archimedean Riesz spaces. Grobler and Labuschagne, on the other hand, in \cite{GL1} and \cite{GL2}, give proofs that avoid representations. In this note we take the opposite approach and show that making more initial use of well-known representations allows a short construction, albeit based on the work involved in proving those classical representation theorems.

Following Troitsky, we call a bilinear map $\phi:E\times F\to G$ \emph{bi-injective} if $\phi(x,y)=0$ implies that either $x=0$ or $y=0$.
We start our construction by pointing out that, starting with Archimedean Riesz spaces $E$ and $F$, it is simple, using a suitable representation theorem, to find an Archimedean Riesz space $G$ and a bi-injective Riesz bimorphism $\phi:E\times F\to G$. Using the classical Krein-Kakutani representation for order unit spaces, we show that the Riesz subspace of $G$ generated by $\phi(E\times F)$ is, essentially, independent of the choice of $G$ and $\phi$. The same techniques allow easy deduction of the basic properties of the tensor product.

The author would like to thank Vladimir Troitsky for suggesting this approach to the construction of the tensor product of Archimedean Riesz spaces as well as for pointing out minor errors in the first draft of this paper.

\section{The starting point}
Recall that if $E, F$ and $G$ are Archimedean Riesz spaces then a linear operator $T:E\to G$ is a \emph{Riesz homomorphism} if $T(x\vee y)=Tx\vee Ty$ for all $x,y\in E$. The bilinear map $\phi:E\times F\to G$ is a \emph{Riesz bimophism} if $\phi(\cdot,y):E\to G$ and $\phi(x,\cdot):F\to G$  are Riesz homomorphisms for all $y\in F_+$ and $x\in E_+$ respectively. We denote the constantly one (resp. zero) function on a set $X$ by $\1_X$ (resp. $\0_X$), and omit the subscript if the domain is unambiguous. If $X$ and $Y$ are topological spaces then the \emph{natural bilinear embedding}  $\chi:C(X)\times C(Y)\to C(X\times Y)$ satisfies $\chi(f,g)(x,y)=f(x)g(y)$ for $x\in X, y\in Y, f\in C(X)$ and $g\in C(Y)$.

If $\Sigma$ is a topological space, let $S(\Sigma)$ be the set of a continuous real-valued functions $f$ defined on dense open subsets $D(f)$ of $\Sigma$. After quotienting out the relation $f\sim g$ if $f$ and $g$ agree on $D(f)\cap D(g)$, $S(\Sigma)$ becomes an Archimedean Riesz space under the usual pointwise linear and order definitions on the appropriate domains. If $\Sigma_1$ and $\Sigma_2$ are two topological spaces, $f_1\in S(\Sigma_1)$ and $f_2\in S(\Sigma_2)$ then $D(f_1)\times D(f_2)$ is a dense open subset of $\Sigma_1\times \Sigma_2$ and $\chi$ may be extended by defining $\chi(f_1,f_2)(\sigma_1,\sigma_2)=f_1(\sigma_1)f_2(\sigma_2)$ on $D(f_1)\times D(f_2)$. It is clear that $\chi:S(\Sigma_1)\times S(\Sigma_2)\to S(\Sigma_1\times \Sigma_2)$ is a bi-injective Riesz bimorphism. As $S(\Sigma_1\times \Sigma_2)$ is an Archimedean Riesz space we may form the Riesz subspace generated by $S(\Sigma_1)\otimes S(\Sigma_2)$.

 In \cite{BW} we pointed out that if one takes the Ogasawara-Maeda representation $x\mapsto x^\wedge$, \S50 of \cite{LZ}, of an Archimedean Riesz space $E$ in $C^\infty(\Sigma_1)$, where $\Sigma_1$ is a Stonean space and compose this with the restriction map of $C^\infty(\Sigma_1)$ into $S(\Sigma_1)$ which simply restricts $f$ to the dense open set on which it is real-valued, then we have a representation of $E$ in $S(\Sigma_1)$.\footnote{It is possible to give a direct and much simpler construction of a representation of an Archimedean Riesz space in a space $S(\Sigma)$, which will be the subject of another paper.} If we similarly have a representation of $F$ in $S(\Sigma_2)$ then we may form the Riesz subspace generated by $\chi(E^\wedge,F^\wedge)$ in $S(\Sigma_1\times \Sigma_2)$. Because of the uniqueness that Fremlin proved, this gives us Fremlin's tensor product, but in \cite{BW} we were at pains to point out that we had not given a construction of the tensor product. This concrete description of the tensor product proved useful not only in \cite{BW}, but also in \cite{Z}. The current paper shows how this observation can actually be used as a  starting point for the construction of the tensor product.

 \section{Uniqueness and properties}

 Our construction will start by considering spaces with a strong order unit, The following result appears in \cite{S}, Proposition III.2.1 but we include a proof for completeness.

\begin{theorem}\label{Schaefer}
Let $X, Y$ and $Z$ be compact Hausdorff spaces, and let $E\subseteq C(X)$ and $F\subseteq C(Y)$  be supremum norm dense Riesz subspaces containing the constants. Let $\phi:E\times F\to C(Z)$ be a Riesz bimorphism then it can be factored as
\[E\times F \xrightarrow{\chi} C(X\times Y)\xrightarrow{T} C(Z),\]
where $\chi$ is the natural bilinear embedding and $T$ is a lattice homomorphism which is uniquely determined by $\phi$.

\begin{proof}
By continuity we may assume that $E=C(X)$ and $F=C(Y)$. If $z\in Z$ then the maps $\phi(\cdot,\1_Y)(z):C(X)\to\Real$ and $\phi(\1_X,\cdot):C(Y)\to \Real$ are lattice homomorphisms, so there are $\lambda_z\in\Real_+$ and $x_z\in X$ with $\phi(f,\1_Y)(z)=\lambda f(x_z)$ for all $f\in C(X)$. Similarly there are $\mu_z\in \Real$ and $y_z\in Y$  with $\phi(\1_X,g)(z)=\mu g(y_z)$ for all $g\in C(Y)$. Considering $\phi(\1_X,\1_Y)$ we see that $\lambda_z=\mu_z$. We claim that $\phi(f,g)(z)=\lambda_z f(x_z)g(y)z)$ for all $f\in C(X), g\in C(Y)$. First  note that if $g(y_z)=0$ then $\phi(f,g)(z)=0$ as
\[|\phi(f,g)|(z)\le \phi(|f|,|g|)(z)\le \phi(\|f\|_\infty \1_X,|g|)(z)=\|f\|_\infty \lambda_z |g|(y_x)=0.\]
Hence 
\begin{align*}
\phi(f,g)(z)&=\phi(f,g(y_z)\1_Y)(z)+\phi(f,g-g(y_z)\1_Y)(z)\\
&=g(y_z)\lambda_z f(x_z)+0=\lambda_z f(x_z)g(y_z).
\end{align*}

Define $T$ on $\chi(C(X),C(Y))$ by $T(\chi(f,g))=\phi(f,g)$. Extend $T$ linearly to the linear span of $\chi(C(X),C(Y))$, i.e. $C(X)\otimes C(X)$, and then by continuity to the whole of $C(X\times Y)$, noting that $C(X)\otimes C(Y)$ is dense in $C(X\times Y)$ for the supremum norm. We claim that if $h\in C(X\times Y)$ then $Th(z)=\lambda_z h(x_z,y_z)$. This is immediate from the definition if $h=\chi(f,g)$ then use linearity followed by continuity. It is now evident that $T$ is a lattice homomorphism. The uniqueness is clear.
\end{proof}

\end{theorem}

\begin{lemma}\label{Urysohn}
Suppose that $X$ is a compact Hausdorff space, $E$ a supremum norm dense Riesz subspace of $C(X)$, containing the constants. If $U$ is a non-empty open set in $X$ containing $w\in X$ then there is $h\in E_+$ with $\0\le h\le \1$, $h(\omega)=1$ and $h_{|\Omega\setminus U}\equiv 0$.
\begin{proof}
By Urysohn's lemma there is $f\in C(\Omega)$ with $-\frac12\1\le f\le \frac32\1$, $f(\omega)=\frac32$ and $f$ constantly $-\frac12$ on $\Omega\setminus U$. As $E$ is dense in $C(\Omega)$ for the supremum norm, there is $g\in E$ with $\|g-f\|_\infty\le \frac12$. This means that $g(\sigma)\le 0$ for $\sigma\in\Omega\setminus U$, and $g(\omega)\ge 1$. Now take $h=g^+\wedge \1$.
\end{proof}
\end{lemma}

\begin{proposition}With the notation of Theorem \ref{Schaefer}, if we also assume that $\phi(\1_X, \1_Y)=\1_Z$ and that $\phi$ is bi-injective then $T$ is injective.
\begin{proof}
As $\chi(\1_X,\1_Y)=\1_{X\times Y}$, $T(\1_{X\times Y})=\1_Z$, so there is a continuous map $\pi:Z\to X\times Y$ with $Tf=f\circ \pi$ for $f\in C(X\times Y)$, \cite{MN} Theorem 3.2.12. We claim that $\pi$ is onto, from which $T$ being injective is immediate. If not, let $(x,y)\in (X\times Y)\setminus \pi(Z)$.  As $Z$ is compact and $\pi$ continuous, this complement is open in $X\times Y$. Thus there are open sets $U\subset X$ and $V\subset Y$ with $x\in U$, $y\in V$ and $(U\times V)\cap \pi(Z)=\emptyset$. Use Lemma \ref{Urysohn} to construct $f\in E$ with $f(x)=1$ and $f_{|X\setminus U}\equiv 0$ and also $g\in F$ with $g(y)=1$ and $g_{|Y\setminus V}\equiv 0$. Now $\chi(f,g)$ vanishes on $\pi(Z)$, so that $\phi(f,g)=T(\chi(f,g))=0$ even though neither $f$ nor $g$ is zero, contradicting $\phi$ being bi-injective.
\end{proof}
\end{proposition}

If $E$ and $F$ are Archimedean Riesz spaces and $\phi$ is a bi-injective Riesz bimorphism of $E\times F$ into an Archimedean Riesz space $G$, then we call the Riesz subspace of $G$ generated by $\phi(E,F)$ \emph{a tensor product} of $E$ and $F$ and denote it by $E\ftp_\phi F$. At this stage we make no assertions about its categorical properties, uniqueness or even its existence. Recall, \cite{J}, Theorem 2.2.11, that the Riesz subspace generated by a linear subspace $A$ of a Riesz space is $\vee(\wedge A)=\wedge(\vee A)$, where $\wedge A$ (resp. $\vee A$) is the set of all finite suprema (resp. infima) from $A$.

\begin{theorem}\label{oucase} Let $X$ and $Y$ be compact Hausdorff spaces, $E\subseteq C(X)$ and $F\subseteq C(Y)$ be dense Riesz subspaces containing the respective constants. Let also $G$ be an Archimedean Riesz space and $\phi:E\times F\to G$ be a Riesz bimorphism. Let also $\chi:E\times F\to C(X\times Y)$ be the natural bilinear embedding. Then there is a unique Riesz homomorphism $T:E\ftp_\chi F\to G$, such that $\phi=T\circ \chi$. Furthermore, if $\phi$ is bi-injective then $T$ is a Riesz isomorphism of $E\ftp_\chi F$ onto $E\ftp_\phi F$.
\begin{proof}
Let $u=\phi(\1_X,\1_Y)\in  G$ and let $Z$ be a compact Hausdorff space such that the principal ideal $G_u$ is Riesz isomorphic to a dense Riesz subspace $H$ of $C(Z)$, with $u$ corresponding to $\1_Z$. We know that there is a lattice homomorphism $T:C(X\times Y)\to C(Z)$ with $\phi=T\circ\chi$, which is injective if $\phi$ is bi-injective. $T$ maps $E\ftp_\chi F$ into $H$, because $T$ is a Riesz homomorphism, $T(\chi(E,F))=\phi(E,F)\subset H$, the characterization of a generated Riesz subspace and $H$ being a Riesz subspace of $C(Z)$. If $\phi$ is bi-injective then  $T_{|E\ftp_\chi F}$ is an injective lattice homomorphism of $E\ftp_\chi F$ onto the Riesz subspace generated by $\phi(E,F)$, $E\ftp_\phi F$. The uniqueness follows from the uniqueness in Theorem \ref{Schaefer}.
\end{proof}
\end{theorem}

\begin{corollary}\label{ouunique} Let $X$ and $Y$ be compact Hausdorff spaces, $E\subseteq C(X)$ and $F\subseteq C(Y)$ be dense Riesz subspaces containing the respective constants. Let also $G_1, G_2$ be Archimedean Riesz spaces and $\phi_i:E\times F\to G_i$ be a bi-injective Riesz bimorphisms then there is a unique Riesz isomorphism $T:E\ftp_{\phi_1} F\to E\ftp_{\phi_2} F$ such that $T\circ\phi_1=\phi_2$.
\begin{proof}
Go via $E\ftp_\chi F$!
\end{proof}
\end{corollary}

Before the next proof note that if $E$ is an Archimedean Riesz space, $A$ a positively generated linear subspace and $G$ the Riesz subspace that it generates then $A$ is cofinal in $G$. This applies in particular to $E\otimes_\phi F$ inside $E\ftp_\phi F$. If $E$ and $F$ have order units then the following result reduces, modulo the Krein-Kakutani theorem, to Corollary \ref{ouunique}.

\begin{theorem}\label{uniqueness}
If $E, F, G_1$ and $G_2$ are Archimedean Riesz spaces and $\phi_i:E\times F\to G_i$ are bi-injective Riesz bimorphisms, then the Riesz subspaces of $G_i$ generated by $\phi_i(E,F)$ are isomorphic via a unique lattice isomorphism $T$ with $T\circ \phi_1=\phi_2$.
\begin{proof}If $u\in E_+$ and $v\in F_+$ then we know there is a unique lattice isomorphism $T_{u,v}:E_u\ftp_{\phi_1} F_v\to E_u\ftp _{\phi_2} F_v$ with $T_{u,v}\circ \phi_1=\phi_2$. If $u\le u'\in E_+$ and $v\le v'\in F_+$ then we may also find $T_{u',v'}:E_{u'}\ftp_{\phi_1} F_{v'}\to E_{u'}\ftp_{\phi_2} F_{v'}$ with $T_{u',v'}\circ\phi_1=\phi_2$. By the uniqueness, $T_{u',v'|E_u\ftp_{\phi_1}F_v}=T_{u,v}$. We may therefore unambiguously define $T:E\ftp_{\phi_1} F\to E\ftp_{\phi_2}$ to extend all possible $T_{u,v}$. It is then routine to check that this makes $T$ be  a linear lattice isomorphism.
\end{proof}
\end{theorem}

The preceding proof could have been formulated as an inverse limit argument as used by Fremlin and Schaefer, but we prefer this phrasing  of the argument.
The proof of the following lemma was given in \cite{BW} and also outlined in  \S2.

\begin{lemma}\label{startpoint}For all Archimedean Riesz spaces $E$ and $F$ there is an Archimedean Riesz space $G$ and a bi-injective Riesz bimorphism $\phi:E\times F\to G$.
\end{lemma}

\begin{definition}If $E$ and $F$ are Archimedean Riesz spaces then \emph{the tensor product}, $E\ftp F$, is the Riesz space generated by $\phi(E\times F)$ where $G$ is an Archimedean Riesz space and $\phi:E\times F\to G$ is a bi-injective Riesz bimorphism.
\end{definition}

By Lemma \ref{startpoint} there is always at least one choice of $E\ftp F$ and by Theorem \ref{uniqueness} all such choices are isomorphic in a manner that fixes the (copies of the) algebraic tensor product $E\otimes F$. It is natural now to regard $E\otimes F$ as a subspace of $E\ftp F$. From now on we denote the natural embedding of $E\times F$ into $E\otimes F\subset E\ftp F$ by $\otimes$. We must now start on the task of obtaining the properties of $E\ftp F$.

With our definition the following Corollary 4.5 of \cite{F} is obvious.

\begin{proposition}If $E$ and $F$ are Archimedean Riesz spaces with Riesz subspaces $E_0$ and $F_0$ respectively, then $E_0\ftp F_0$ can be identified with the Riesz subspace of $E\ftp F$ generated by $E_0\otimes F_0$.
\end{proposition}

\begin{theorem}Let $E, F$ and $G$ be Archimedean Riesz spaces and $\phi:E\times F\to G$ is a Riesz bimorphism then there is a unique Riesz homomorphism $T:E\ftp F\to G$ such that $\phi=T\circ\otimes$.
\begin{proof}
For $u\in E_+$ and $v\in F_+$ we know from Theorem \ref{oucase} that there is a unique Riesz homomorphism $T_{u,v}:E_u\ftp F_v\to G_{\phi(u,v)}\subseteq G$ with $\phi=T\circ\otimes$ on the domain. As in the proof of Theorem \ref{uniqueness} we may consistently piece together these $T_{u,v}$ to give a Riesz homomorphism $T:E\ftp F\to G$ with $\phi=T\circ\otimes$.
\end{proof}
\end{theorem}

Now for some density results.

\begin{proposition}
Let $E$ and $F$ be Archimedean Riesz spaces. If $a\in E\ftp F$ there are $u\in E_+$ and $v\in F_+$ such that for all $\epsilon>0$ there is $b\in E\otimes F$ such that $|a-b|\le \epsilon u\otimes v$.
\begin{proof}
Pick $u$ and $v$ such that $w\in E_u\ftp F_v$ and utilize the uniqueness that we proved above. Represent $E_u$ in $C(X)$ and $F_v$ in $C(Y)$, where $X$ and $Y$ are compact Hausdorff spaces, via Riesz isomorphisms $\pi:E_u\to  C(X)$ and $\rho:F_v\to  C(Y)$ with $\pi(u)=\1_X$ and $\rho(v)=\1_Y$ and with $\pi(E_u)$ dense in $C(X)$ and $\rho(F_v)$ dense in $C(Y)$. Now $\pi(E_u)\otimes \rho(F_v)$ is dense in $C(X)\otimes C(Y)$ for the supremum norm. In its turn, $(C(X)\otimes C(Y)$ is dense in $C(X\times Y)$, being a unital separating subalgebra, so $\pi(E_u)\otimes \rho(F_v)$ is dense in $C(X\times Y)$. Using Theorem \ref{oucase} and noting that $\pi(u)\otimes \rho(v)=\1_{X\times Y}$ gives the result.
\end{proof}
\end{proposition}

The next result is slightly more powerful than Fremlin's result in Theorem 4.2 (iv) of \cite{F}.

\begin{proposition}
Let $E$ and $F$ be Archimedean Riesz spaces. If $0\ne c\in (E\ftp F)_+$ then $c=\sup\{a\otimes b: a\in E_+, b\in F_+, a\otimes b\le c\}$
\begin{proof}
As $E\otimes F$ is cofinal in $E\ftp F$ it is not difficult to find $u\in E_+$ and $v\in F_+$ such that $c\le u\otimes v$. We may now restrict our attention to $E_u\ftp F_v$. Identify $E_u$ and $F_v$ with dense Riesz subspaces of $C(X)$ and $C(Y)$ respectively, which contain the constants. At each point $(x,y)\in X\times Y$ with $c(x,y)>0$, we may use continuity to find, for each $n\in\Natural$,  open sets $U_n\subset X$ and $V_n\subset Y$ with $c_{|U\times V}\ge (1-\frac1{n+1})c(x,y)>0$. Now use Lemma \ref{Urysohn} to find $a_n\in (E_u)_+$ and $b_n\in (F_v)_+$ with $a_n$ vanishing on $X\setminus U$, $a_n\le \1_X$, $a_n(x)=1$,  $b_n$ vanishing on $Y\setminus V$, $b_n\le (1-\frac1{n+1})c(x,y)\1_Y$ and $b_n(y)=(1-\frac1{n+1})c(x,y)$. This forces $0\le a_n\otimes b_n\le c$ and $a_n(x)n_n)(y)=(1-\frac1{n+1})c(x,y)$. Hence $\sup\{a(x)b(y):a\in E_+, b\in F_+, a\otimes b \le c\}=c(x,y)$ whenever $c(x,y)>0$. The result is now immediate.
\end{proof}
\end{proposition}


\begin{bibdiv}
\begin{biblist}[\resetbiblist{9}]
\bib{BW}{article}{
   author={Buskes, G. J. H. M.},
   author={Wickstead, A. W.},
   title={Tensor products of $f$-algebras},
   journal={Mediterr. J. Math.},
   volume={14},
   date={2017},
   number={2},
   pages={Paper No. 63, 10},
   issn={1660-5446},
   review={\MR{3619425}},
   doi={10.1007/s00009-017-0841-x},
}

\bib{F}{article}{
   author={Fremlin, D. H.},
   title={Tensor products of Archimedean vector lattices},
   journal={Amer. J. Math.},
   volume={94},
   date={1972},
   pages={777--798},
   issn={0002-9327},
   review={\MR{0312203}},
   doi={10.2307/2373758},
}
\bib{GL1}{article}{
   author={Grobler, J. J.},
   author={Labuschagne, C. C. A.},
   title={The tensor product of Archimedean ordered vector spaces},
   journal={Math. Proc. Cambridge Philos. Soc.},
   volume={104},
   date={1988},
   number={2},
   pages={331--345},
   issn={0305-0041},
   review={\MR{0948918}},
   doi={10.1017/S0305004100065506},
}

\bib{GL2}{article}{
   author={Grobler, J. J.},
   author={Labuschagne, C. C. A.},
   title={An $f$-algebra approach to the Riesz tensor product of Archimedean
   Riesz spaces},
   journal={Quaestiones Math.},
   volume={12},
   date={1989},
   number={4},
   pages={425--438},
   issn={0379-9468},
   review={\MR{1021941}},
}


\bib{J}{book}{
   author={Jameson, Graham},
   title={Ordered linear spaces},
   series={Lecture Notes in Mathematics},
   volume={Vol. 141},
   publisher={Springer-Verlag, Berlin-New York},
   date={1970},
   pages={xv+194},
   review={\MR{0438077}},
}


\bib{LZ}{book}{
   author={Luxemburg, W. A. J.},
   author={Zaanen, A. C.},
   title={Riesz spaces. Vol. I},
   series={North-Holland Mathematical Library},
   publisher={North-Holland Publishing Co., Amsterdam-London; American
   Elsevier Publishing Co., Inc., New York},
   date={1971},
   pages={xi+514},
   review={\MR{0511676}},
}
\bib{MN}{book}{
   author={Meyer-Nieberg, Peter},
   title={Banach lattices},
   series={Universitext},
   publisher={Springer-Verlag, Berlin},
   date={1991},
   pages={xvi+395},
   isbn={3-540-54201-9},
   review={\MR{1128093}},
   doi={10.1007/978-3-642-76724-1},
}
\bib{S}{article}{
   author={Schaefer, H. H.},
   title={Aspects of Banach lattices},
   conference={
      title={Studies in functional analysis},
   },
   book={
      series={MAA Stud. Math.},
      volume={21},
      publisher={Math. Assoc. America, Washington, DC},
   },
   isbn={0-88385-121-0},
   date={1980},
   pages={158--221},
   review={\MR{0589416}},
}


\bib{Z}{article}{
   author={Zabeti, O.},
   title={Fremlin tensor product respects the unbounded convergences},
   journal={(submitted)},
  }\end{biblist}
\end{bibdiv}
\end{document}